\PassOptionsToPackage{dvipsnames}{xcolor}

\documentclass{ifacconf}

\makeatletter
\let\old@ssect\@ssect 
\makeatother

\usepackage[pdftex]{graphicx}
    \graphicspath{{Figures/}}                   
    \DeclareGraphicsExtensions{.pdf,.jpg,.png}  
    \usepackage{svg}                            
    \usepackage{float}                          

\usepackage[dvipsnames]{xcolor}                 
\usepackage{tikz}                               
    \usetikzlibrary{calc,patterns,angles,quotes}

\usepackage{amsmath}                            
\usepackage{amssymb}                            
\usepackage{mathtools}                          
\usepackage{breqn}                              
\usepackage{array}                              
\usepackage{nicefrac}
\usepackage{blkarray, bigstrut}                 

\usepackage[version=4]{mhchem}                  

\usepackage{dutchcal}

\usepackage{natbib}                             
\usepackage{siunitx}                            
\usepackage{textcomp}                           

\usepackage{color}                              

\usepackage{hyperref}                           
    
\newcommand{\R}{\mathbb{R}}                     
\newcommand{\doubleH}{\mathbb{H}}               
\newcommand{\mb}[1]{\boldsymbol{#1}}            
\newcommand{\mcal}{\mathcal}                    
\newcommand{\Id}{\mathbb{I}}
\newcommand{\quat}[1]{\mathfrak{#1}}            
\newcommand{\Bmat}[1]{\boldsymbol{#1}}          

\newcommand{\supfirst}{$1^{\mathrm{st}}$}       
\newcommand{\supsecond}{$2^{\mathrm{nd}}$}      

\DeclareMathOperator{\diag}{diag}
\DeclarePairedDelimiter\abs{\lvert}{\rvert}
\DeclarePairedDelimiter\norm{\lVert}{\rVert}

\makeatletter
\def\@ssect#1#2#3#4#5#6{%
  \NR@gettitle{#6}
  \old@ssect{#1}{#2}{#3}{#4}{#5}{#6}
}

\begin{document}
\begin{frontmatter}

\title{Attitude Trajectory Optimization \\ and Momentum Conservation \\ with Control Moment Gyroscopes\thanksref{footnoteinfo}}  

\thanks[footnoteinfo]{This material is based upon work supported by the National Science Foundation Graduate Research Fellowship Program under Grant No. DGE 1650115. Any opinions, findings, and conclusions or recommendations expressed in this material are those of the author(s) and do not necessarily reflect the views of the National Science Foundation.}
\thanks[CorAuthorInfo]{Corresponding Author (email: Thomas.Dearing@colorado.edu)}

\author[First]{Thomas~L.~Dearing\thanksref{CorAuthorInfo}} 
\author[First]{John~Hauser}
\author[Second]{Christopher~Petersen}
\author[First]{Marco~M.~Nicotra}
\author[First]{Xudong~Chen}

\address[First]{University of Colorado Boulder, ECEE Department,
   Boulder, CO 80309 USA (e-mails: Firstname.Lastname@colorado.edu).}
\address[Second]{University of Florida, Department of Mechanical and Aerospace Engineering,
   Gainesville, FL 32611 USA (e-mail: c.petersen1@ufl.edu)}

\begin{abstract}
In this work, we develop a numerically tractable trajectory optimization problem for rest-to-rest attitude transfers with CMG-driven spacecraft. First, we adapt a specialized dynamical model which avoids many of the numerical challenges (singularities) introduced by common dynamical approximations. To formulate and solve our specialized trajectory optimization problem, we design a locally stabilizing Linear Quadratic (LQ) regulator on the system's configuration manifold then lift it into the ambient state space to produce suitable terminal and running LQ cost functionals. Finally, we examine the performance benefits and drawbacks of solutions to this optimization problem via the PRONTO solver and find significant improvements in maneuver time, terminal state accuracy, and total control effort. This analysis also highlights a critical shortcoming for objective functions which penalize only the norm of the control input rather than electrical power usage.
\end{abstract}



\end{frontmatter}



\section{Introduction}
A fundamental challenge when designing a spacecraft is achieving a balance between the available onboard power and fuel and the performance capabilities necessary to complete the mission (\cite{Larson1999}). Naturally, the efficiency of the spacecraft's Attitude Determination and Control System (ADCS) is paramount in this balance, as each improvement in efficiency enables more resources dedicated to the mission objective (additional sensors, faster processors, etc.). Naturally, \emph{optimal} control strategies prove invaluable in this context as the objective function can be tailored to the specific mission. 

While optimal control strategies have been identified for simple thruster-driven spacecraft, platforms using more efficient, prevalent, and complex Control Moment Gyroscopes (CMG's) for attitude control have presented significant challenges for conventional optimization approaches:
\begin{enumerate}
    \item Dynamics evolving on a non-Euclidean manifold
    \item High state depth, integrator order, and nonlinearity
    \item Numeric challenges from prevalent approximations
\end{enumerate}
In particular, Momentum-Exchange Devices (MED's) like CMG's operate via the conservation of total angular momentum: a nonlinear constraint which shapes the system's unique state manifold. To accommodate the sophisticated dynamics induced by MED's, existing optimal control solutions predominantly employ \emph{approximated} dynamics to reduce numerical complexity. For example, the varied optimization approaches presented in 
\cite{Lee2017} (Indirect Single Shooting), 
\cite{Banerjee2019} (Pseudospectral methods), and
\cite{Wang2020} (Differential Evolution methods) 
all plan maneuvers using a thruster-driven dynamical model (omitting momentum conservation entirely). However, the challenges involved in obtaining and integrating even these solutions with conventional CMG arrays has limited widespread adoption.

In this work, we extend the trajectory optimization approach in \cite{Dearing2021} to a CMG-driven satellite model. In contrast to other approaches, we use a specialized dynamical model that preserves the system's (conserved) momentum exchange physics using the original CMG motor torques as control inputs. Finding this formulation to be substantially more numerically efficient than conventional approximated models for our solver, we then examine optimal solutions for the popular rooftop and pyramid CMG array geometries.


The remaining sections are organized as follows: %
    Section II introduces the notations and physics for MED driven spacecraft, %
    Section III presents our dynamical model and its benefits, %
    Section IV formulates our optimization problem, %
    Section V presents our optimal trajectories and their comparative performance, %
    and Section VI summarizes the paper with concluding remarks.



\section{Satellites and Momentum Exchange}

\subsection{Attitude Representations}
In this work, satellite rotations are modelled using quaternions both for their computational efficiency and lack of coordinate singularities. A \emph{quaternion} $\quat{q} \in \doubleH$ is a hypercomplex number of the form
$$
\quat{q} \coloneqq q_s + \underbrace{q_x \mb{i} + q_y \mb{j} + q_z \mb{k}}_{ \quat{q}_v \coloneqq } \, ,
$$
with real (scalar) part $q_s \coloneqq \mathrm{Re}(\quat{q}) \in \R$ and imaginary (vector) part $\quat{q}_v \coloneqq \mathrm{Im}(\quat{q})$ written using the complex basis $\mb{i},\mb{j}$, and $\mb{k}$. For simplicity, $\quat{q}$ is often written as the vector $\mb{q} \coloneqq [q_s;\mb{q}_v] \in \R^4$ with $\mb{q}_v \coloneqq [q_x;q_y;q_z] \in \R^3$, where we use the notation $[\mb{a};\mb{b}]$ to denote vertically concatenated vectors. Correspondingly, the quaternion product ``$\circ$'' admits the following vector notation equivalent:
\begin{align}
    \quat{h} &= \quat{q} \circ \quat{p}, \qquad\qquad\qquad\qquad\qquad \quat{h},\quat{p},\quat{q} \in \doubleH, \notag\\
    \begin{bmatrix}
        h_s \\
        \mb{h}_v
    \end{bmatrix}
      &= 
    \begin{bmatrix}
        q_s p_s - \mb{q}_v^\top \mb{p}_v \\
        q_s \mb{p}_v + p_s \mb{q}_v + \mb{q}_v \times \mb{p}_v
    \end{bmatrix}, \label{eq_quat_prod_def}\\
      &= 
    \underbrace{
    \begin{bmatrix}
             q_s &       -\mb{q}_v^\top      \\
        \mb{q}_v & q_s \Id_3 + \widehat{\mb{q}}_v
    \end{bmatrix}
    }_{O_L(\mb{q}) \coloneqq}
    \begin{bmatrix}
        p_s \\
        \mb{p}_v
    \end{bmatrix}
    = 
    \underbrace{
    \begin{bmatrix}
        p_s &       -\mb{p}_v^\top      \\
        \mb{p}_v & p_s \Id_3 - \widehat{\mb{p}}_v
    \end{bmatrix}
    }_{O_R(\mb{p}) \coloneqq}
    \begin{bmatrix}
        q_s \\
        \mb{q}_v
    \end{bmatrix}, \notag
\end{align}
where $\Id_3$ denotes the $3\times3$ identity matrix and $O_L(\mb{q})$ and $O_R(\mb{p})$ are the (orthogonal) matrix representations of the quaternion product from the left (by $\quat{q}$) and from the right (by $\quat{p}$), respectively. Additionally, the hat operator
$$
\widehat{\mb{\omega}} \coloneqq
\begin{bmatrix}
    \phantom{-}0        & -\omega_z           & \phantom{-}\omega_y \\
    \phantom{-}\omega_z & \phantom{-}0        & -\omega_x \\
    -\omega_y           & \phantom{-}\omega_x & \phantom{-}0
\end{bmatrix},
$$
yields the matrix representation of the cross product $\widehat{\mb{\omega}} \mb{v} = \mb{\omega} \times \mb{v}$. Finally, the \emph{conjugate} for a quaternion $\quat{q}$ is defined as $\quat{q}^* \coloneqq q_s - \quat{q}_v$ and the quaternion norm is given by $\norm{\quat{q}}^2 \coloneqq \quat{q}^* \circ \quat{q}$ and agrees with the usual Euclidian norm on $\R^4$ (e.g. $\norm{\quat{q}} \equiv \norm{\mb{q}}$). We refer the reader to \cite{Schwab2002} and \cite{DeRuiter2013} for further details.

Like conventional rotation matrices, unit quaternions ($\norm{\quat{q}} = 1$) can be used to represent spacecraft attitudes and rotations. Specifically, a vector $\mb{v}^b\in \R^3$ written in the satellite's body frame $\mcal{F}_b$ can be transformed to the space-fixed inertial frame $\mcal{F}_i$ using either the quaternion $\quat{q} \in \doubleH$ or the rotation matrix $C(\mb{q}) \in \mathrm{SO}(3)$ as follows:
\begin{subequations}
\label{eq_rot_transforms}
\begin{align}
    \widetilde{\mb{v}}^{\,i} &= \quat{q} \circ \widetilde{\mb{v}}^{\,b} \circ \quat{q}^*, \label{eq_quat_transform}\\
    \mb{v}^i &= C(\mb{q}) \, \mb{v}^b. \label{eq_rotm_transform}
\end{align}
\end{subequations}
where $\widetilde{\mb{v}} \coloneqq [0;\mb{v}]$ and
$$
C(\mb{q}) \coloneqq q_s^2\, \Id_{3} + 2 q_s \widehat{\mb{q}}_v + \mb{q}_v \mb{q}_v^\top + \widehat{\mb{q}}_v^{\phantom{.}2}.
$$

\subsection{Momentum Exchange Devices}
To minimize the use of consumable fuels, satellite attitude is nominally controlled using only renewable electric power. While specialized devices such as magnetic field torquers can generate weak \emph{external} torques on the satellite body, far greater agility can be achieved using simple electric motors to \emph{internally} redistribute the platform's angular momentum. The simplest such device is the Reaction Wheel (RW): an electric motor mounted to the frame of the satellite with a high-inertia rotor (\cite{Larson1999}). When the motor applies torque to rotate the wheel, the resulting reaction torque of the wheel on the motor frame is used to rotate the satellite. An array of such devices mounted on the satellite's principle inertia axes thus produces a reliable attitude control system.

Unfortunately, RW's are quite inefficient for heavier spacecraft. The mechanical shaft power $P = \tau_w \omega_w$ required for a motor to apply a torque $\tau_w$ increases linearly with the wheel speed $\omega_w$.  Thus, RW torque generation is inefficient at high wheel speeds, while friction effects can also make RW's unreliable at low speeds. As a result, RW control systems require active wheel speed regulation to avoid both effects (\cite{Leve2015}).

\begin{figure}
    \centering
    \includegraphics[width=2.4in]{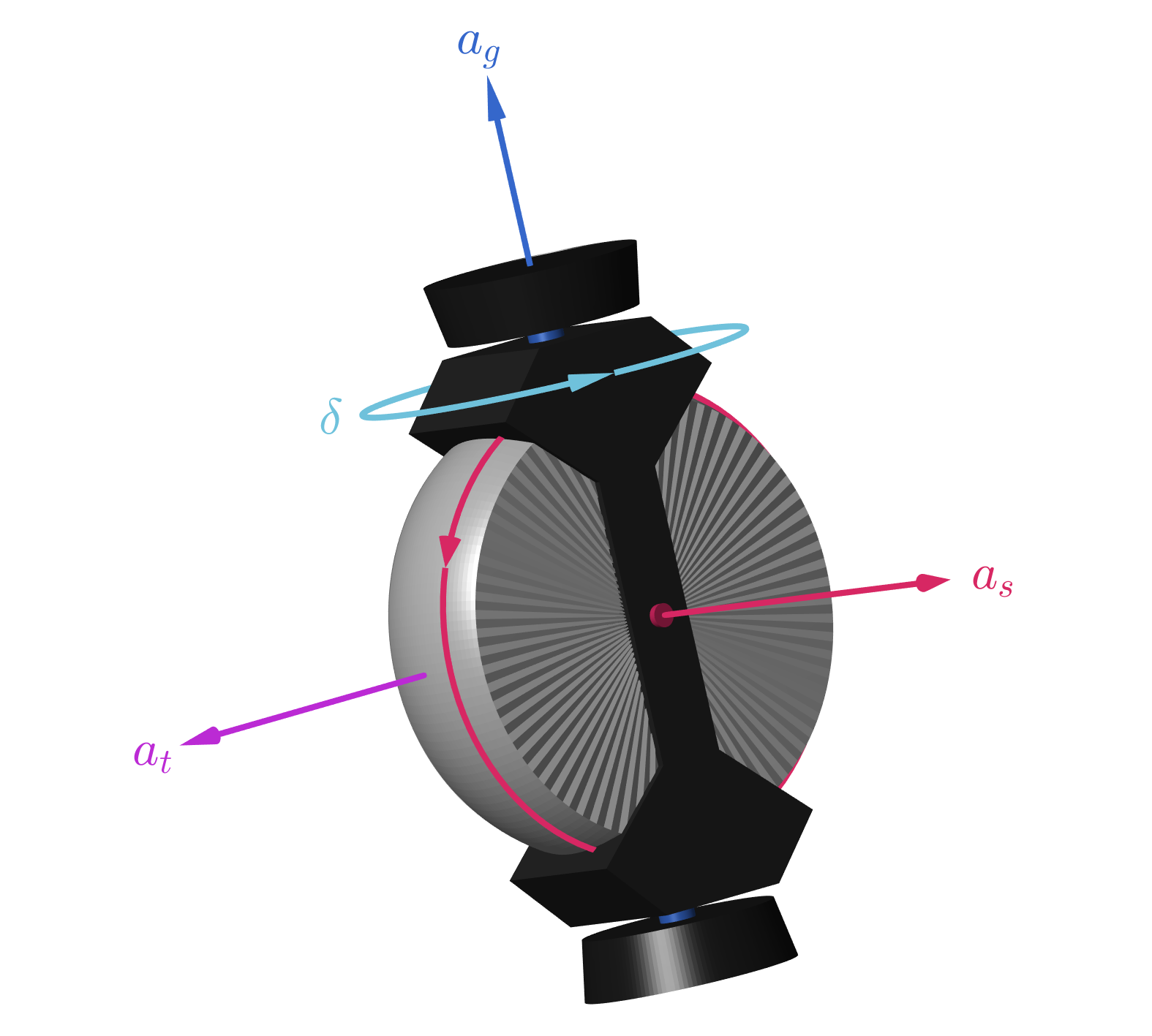}
    \caption{Control Moment Gyroscope coordinate frame.
    \label{Fig_VSCMG}
    }
\end{figure}

Evolving from this design, the Control Moment Gyroscope is a reaction wheel mounted to a rotating gimbal as shown in Fig. \ref{Fig_VSCMG}, where the wheel (red) and gimbal (blue) motors act along the $\mb{a}_s$ and $\mb{a}_g$ axes respectively. Rather than using the motor reaction torques for attitude control, a CMG instead employs the \emph{gyroscopic} reaction torque 
$$
\mb{\tau}_{r} = \dot{\delta} h_w \mb{a}_t,
$$
produced along the \emph{transverse} axis $\mb{a}_t \coloneqq \mb{a}_s \times \mb{a}_g$. Note that $\mb{\tau}_r$ is proportional to the rotation \emph{rate} $\dot{\delta} \in \R$ of the gimbal assembly (\emph{not} the gimbal motor torque $\mb{\tau}_g$) and is amplified by the rotor momentum $h_w \in \R$. This \emph{torque amplification} allows CMG's to efficiently generate larger output torques than RW's (e.g. $\geq \SI{1000}{\newton\meter}$ vs. $\leq \SI{1}{\newton\meter}$).

However, while a RW's torque axis remains fixed in the body frame $\mcal{F}_b$, a CMG's output torque axis $\mb{a}_t$ rotates with the gimbal angle $\delta \in [0,2\pi)$. As such, the available output torque from an array of $m$ CMG's varies with the array's \emph{configuration} $\mb{\delta} := [\delta_1;\cdots;\delta_m]$; a drawback that demands more elaborate control strategies. Following the notation in \cite{Ford2000}, the available torque spaces for the CMG gimbal and wheel motors are spanned respectively by the column spaces of the matrices
\begin{equation}
\begin{aligned}
    A_s &\coloneqq [\mb{a}_{s,1},\cdots,\mb{a}_{s,m}], \quad A_t \coloneqq [\mb{a}_{t,1},\cdots,\mb{a}_{t,m}], \\
\end{aligned}
\end{equation}
which vary with the array configuration $\mb{\delta}$ following
\begin{equation}
\begin{aligned}
    A_s(\mb{\delta}) &\coloneqq A_{s0} \diag( \cos (\mb{\delta})) - A_{t0} \diag( \sin (\mb{\delta})),\\
    A_t(\mb{\delta}) &\coloneqq A_{t0} \diag( \cos (\mb{\delta})) + A_{s0} \diag( \sin (\mb{\delta})),\\
\end{aligned}
\end{equation}
where $\sin()$ and $\cos()$ act entry-wise for vector inputs and $A_{s}(\mb{0}) = A_{s0}$ and $A_{t}(\mb{0}) = A_{t0}$ define the default configuration of the array geometry. For completeness, the (fixed) gimbal axes are also collected in the constant matrix $A_g = [\cdots,\mb{a}_{g,i},\cdots]$.  Under this notation, the column space of $A_t$ describes the available gyroscopic torques from the gimbal motors, while that of $A_s$ describes the (RW) reaction torques available from the wheel motors.

\subsection{Momentum and Inertia of a CMG Array}
Next, we examine the momentum exchange physics of a CMG array to determine the satellite's variable Moment of Inertia (MoI) and body-frame angular momentum. First let $J_B \in \R^{3 \times 3}$ be the constant diagonal inertia of the satellite body in $\mcal{F}_b$ omitting the CMG array. Next, let the array have $m$ CMG's with relative positions $\mb{r}_i \in \R^3$ and orientations in $\mcal{F}_b$ given by the matrices $A_g$, $A_s$, and $A_t$ as above. Let each CMG have mass $m_i$ and principle inertia $J_{g,i}$, $J_{s,i}$, and $J_{t,i} \in \R$ along their gimbal, spin, and transverse axes respectively. Collecting these inertias into the $m \times m$ diagonal matrices $\Bmat{J}_s$, $\Bmat{J}_t$, and $\Bmat{J}_g$ (e.g.  $\Bmat{J}_s \coloneqq \mathrm{diag}( [\cdots, J_{s,i}, \cdots])$) and applying the parallel axis theorem, the satellite's total MoI is assembled as follows:
\begin{subequations}
\label{Eq_JstgJ_definition}
\begin{align}
    J_{stg}(\mb{\delta}) &\coloneqq J + A_g \Bmat{J}_g A_g^\top + A_s \Bmat{J}_s A_s^\top + A_t \Bmat{J}_t A_t^\top, \label{Eq_Jstg_definition}\\
    J &= J_B + \sum_{i=1}^{m} m_i \left( \Id_3 \norm{\mb{r}_i}^2 - \mb{r}_i \mb{r}_i^\top \right), \label{Eq_J_definition}
\end{align}
\end{subequations}
where $\Id_3$ is the $3 \times 3$ identity matrix. Note that the first two terms of \eqref{Eq_Jstg_definition} are constant as the CMG gimbal axes and centers of mass are fixed in the body frame.

To determine the satellite's angular momentum, let the body frame $\mcal{F}_b$ have an angular rotation rate $\mb{\omega} \in \R^3$ (measured in $\mcal{F}_b$) with respect to the inertial frame $\mcal{F}_i$. Since the CMG's also rotate within $\mcal{F}_b$, we collect the individual \emph{relative} (to $\mcal{F}_b$) angular momenta around their spin, gimbal, and transverse axes to form the vectors $\mb{h}_{sr}$, $\mb{h}_{gr}$, and $\mb{h}_{tr} \in \R^m$ respectively (though $\mb{h}_{tr} \equiv \mb{0}$ as the CMG cannot rotate around $\mb{a}_t$ in $\mcal{F}_b$). Adding these together, the satellite's total angular momentum in $\mcal{F}_b$ is
\begin{equation}
\label{Eq_h_constructive_definition}
    \mb{h} = J_{stg} \mb{\omega} + A_s \mb{h}_{sr} + A_g \mb{h}_{gr}.
\end{equation}
Building upon this notation, we may alternatively consider the dual \emph{absolute} CMG momenta $\mb{h}_{ga}$ and $\mb{h}_{sa}$ obtained by including the CMG's angular momentum with respect to $\mcal{F}_i$ (nominally embedded in $J_{stg} \mb{\omega}$). These are given by
\begin{subequations}
\begin{align}
    \mb{h}_{ga}  = \mb{h}_{gr} + \Bmat{J}_g A_g^\top \mb{\omega}, \qquad \mb{h}_{sa} &= \mb{h}_{swa} + \mb{h}_{sga}, \label{Eq_hgadef}\\
    \mb{h}_{swa} &= \mb{h}_{swr} + \Bmat{J}_{sw} A_s^\top \mb{\omega}, \label{Eq_hswadef}\\
    \mb{h}_{sga} &= \mb{h}_{sgr} + \Bmat{J}_{sg} A_s^\top \mb{\omega}, \label{Eq_hsgadef}
\end{align}
\end{subequations}
where the momenta $\mb{h}_{sa}$ and inertia $\Bmat{J}_{s} = \Bmat{J}_{sw} + \Bmat{J}_{sg}$ are partitioned between the CMG gimbal frame and wheel respectively. Noting again that the gimbal frame cannot rotate around $\mb{a}_s$, we have that $\mb{h}_{sgr} \equiv 0$ or, notationally, $\mb{h}_{swr} = \mb{h}_{sr}$. Under these alternate coordinates, \eqref{Eq_JstgJ_definition} and \eqref{Eq_h_constructive_definition} can be rewritten as 
\begin{subequations}
\label{Eq_hJstdef}
\begin{align}
    J_{st}(\mb{\delta}) &\coloneqq J + A_s \Bmat{J}_s A_s^\top + A_t \Bmat{J}_t A_t^\top, \label{Eq_Jstdef} \\
    \mb{h} &\coloneqq J_{st} \mb{\omega} + A_s \mb{h}_{swr} + A_g \mb{h}_{ga}. \label{Eq_hdef}
\end{align}
\end{subequations}
For compactness, we will often use \eqref{Eq_hdef} to convert between $\mb{h}$ and $\mb{\omega}$ via the following transformations: 
\begin{subequations}
\label{Eq_hwbar_definition}
\begin{align}
    \bar{\mb{h}}(\mb{\omega},\mb{\delta},\mb{h}_{swr},\mb{h}_{ga}) &\coloneqq J_{st} \mb{\omega} + A_s \mb{h}_{swr} + A_g \mb{h}_{ga}, \label{Eq_hbar_definition}\\
    \bar{\mb{\omega}}(\mb{h},\mb{\delta},\mb{h}_{swr},\mb{h}_{ga}) &\coloneqq J_{st}^{-1} \left( \mb{h} - A_s \mb{h}_{swr} - A_g \mb{h}_{ga} \right), \label{Eq_wbar_definition}
\end{align}
\end{subequations}
and will often write $\bar{\mb{h}}$ or $\bar{\mb{h}}(\mb{x})$ for \eqref{Eq_hbar_definition}, and $\bar{\mb{\omega}}$ or $\bar{\mb{\omega}}(\mb{x})$ for \eqref{Eq_wbar_definition} respectively. Finally, \eqref{Eq_hbar_definition} can be used to determine the array's $3 \times m$ actuator Jacobian
\begin{align}
    D(\mb{\omega},\mb{\delta},&\mb{h}_{swr}) \coloneqq \frac{\partial \bar{\mb{h}}}{\partial \mb{\delta}} = \frac{\partial J_{st}}{\partial\mb{\delta}} \mb{\omega} + \frac{\partial A_s}{\partial\mb{\delta}} \mb{h}_{swr}, \label{Eq_Ddef} \\
    &= \left[ A_s \diag(A_t^{\top} \mb{\omega}) + A_t \diag(A_s^{\top} \mb{\omega}) \right](\Bmat{J}_t - \Bmat{J}_s), \notag\\
    &\quad - A_t \diag(\mb{h}_{swr}), \notag
\end{align}
frequently used in existing controllers to relate the CMG gimbal rates $\dot{\mb{\delta}}(t)$ to the array's output torque $\mb{\tau}_r$.

\section{Dynamical Model}

\subsection{Attitude Dynamics of a CMG-driven satellite}
With the satellite's momentum exchange fully modelled, we may now adapt the dynamics presented by \cite{Ford2000} to coordinates suitable for our optimization problem. Specifically, our state and control inputs are
\begin{equation}
\begin{aligned}
    \mb{x} &\coloneqq 
    [\, \mb{q} \,; \mb{h}_{swr} \,; \mb{\omega} \,; \mb{\delta} \,; \mb{h}_{ga} \,] \in \R^{3m+7}, \\
    \mb{u} &\coloneqq
    [\, \mb{u}_g  \,; \mb{u}_w \,] \in \R^{2m}.
\end{aligned}
\end{equation}
where $\mb{u}_w, \mb{u}_g \in \R^m$ collect the CMG motor inputs for the wheel and gimbal respectively. In particular, choosing the CMG wheel momenta $\mb{h}_{swr}$ as a state allows wheel speed regulation to be easily incentivized in the cost function. Additionally, choosing the original CMG motor torques as control inputs gives the optimizer maximal control over the array's momentum, allowing the full range of direct and gyroscopic reaction control torques.  Notably, choosing to model the wheel momentum $\mb{h}_{swr}$ as variable (rather than constant) identifies this as a \emph{Variable Speed} CMG (VSCMG) model, though this prevalent distinction is insignificant in practice as VSCMG's and CMG's are mechanically identical. Using the notation $\dot{\mb{\nu}} \coloneqq \mb{f}_{\nu}(\mb{x},\mb{u})$ for $\mb{\nu} \in \{\mb{h},\mb{h}_{swr},\mb{\omega}, \mb{\delta}, \mb{h}_{ga}\}$, the complete body-frame dynamics for a CMG-driven satellite are given by
\begin{subequations}
\label{Eq_Dynamics}
\begin{align}
    \dot{\mb{q}} &=  \nicefrac{1}{2} \, O_L(\mb{q}) \, \widetilde{\mb{\omega}}, \label{Eq_Dynamics_q}\\
    \dot{\mb{h}}_{swr} &=  \Bmat{J}_{sw} \left[  \mathrm{diag}(A_t^\top\mb{\omega}) \mb{f}_{\delta} - A_s^\top \mb{f}_{\omega} \right] + \mb{u}_w, \label{Eq_Dynamics_hwr}\\
    \dot{\mb{\omega}} &=  J_{st,a}^{-1}\left[ \mb{f}_h - D_a \mb{f}_{\delta} - A_g \mb{f}_{hga} - A_s \mb{u}_w \right], \label{Eq_Dynamics_w}\\
    \dot{\mb{\delta}} &=  \Bmat{J}_g^{-1} \mb{h}_g - A_g^{\top} \mb{\omega}, \label{Eq_Dynamics_d}\\
    \dot{\mb{h}}_{ga} &= \mathrm{diag} \left( A_t^{\top} \mb{\omega} \right) \left[ (\Bmat{J}_t - \Bmat{J}_s) A_s^{\top} \mb{\omega} - \mb{h}_{swr} \right] + \mb{u}_g, \label{Eq_Dynamics_hg}
\end{align}
\end{subequations}
where the satellite's angular momentum dynamics are
\begin{equation}
\label{Eq_hdyndef}
    \mb{f}_h(\mb{x},\mb{\tau}_e) \coloneqq \widehat{\bar{\mb{h}}} \mb{\omega} + \mb{\tau}_e,
\end{equation}
and $\mb{\tau}_e \in \R^3$ collects any known external torques on the satellite body (e.g. atmospheric drag). For compactness, \eqref{Eq_Dynamics} uses variations on the satellite's MoI \eqref{Eq_Jstdef} given by
\begin{equation}
\label{Eq_Jstadef}
    J_{st,a}(\mb{\delta}) \coloneqq  J + A_s \Bmat{J}_{sg} A_s^\top + A_t \Bmat{J}_t A_t^\top,
\end{equation}
as well as on the actuator Jacobian \eqref{Eq_Ddef} given by
\begin{align}
    D_a(\mb{\omega},\mb{\delta},&\mb{h}_{swa}) \coloneqq \frac{\partial J_{st,a}}{\partial\mb{\delta}} \mb{\omega} + \frac{\partial A_s}{\partial\mb{\delta}} \mb{h}_{swa}, \label{Eq_Dadef}\\
    &= \left[ A_s \diag(A_t^{\top} \mb{\omega}) + A_t \diag(A_s^{\top} \mb{\omega}) \right](\Bmat{J}_t - \Bmat{J}_{sg}) \notag\\
    &\quad - A_t \diag(\mb{h}_{swa}). \notag
\end{align}
Finally, equations \eqref{Eq_Dynamics} are ordered by computational dependency with each dynamical block requiring results only from lower blocks (e.g. $\mb{f}_{\omega}$ depends on $\mb{f}_{\delta}$ but not $\mb{f}_{hswr}$).

While the high state count and nonlinearity of the dynamics \eqref{Eq_Dynamics} would normally make it unsuitable for optimization methods, current literature reveals that approaches using approximated models have encountered far more substantial challenges. We review these approximations, their drawbacks, and the comparative benefits offered by \eqref{Eq_Dynamics} in the following sections.


\subsection{Dynamical Constraints}
A critical aspect of MED-driven spacecraft models that is easily lost in approximation is the conservation of the spacecraft's total inertial angular momentum. This constraint, along with the geometric restriction of the attitude $\mb{q}(t)$ to the unit sphere $\mathrm{S}^3$, are given by
\begin{subequations}
\label{Eq_PhysConstraints}
\begin{align}
    1 &= \norm{\mb{q}}, \label{Eq_PhysConstraints_q} \\
    \mb{h}_0 &= C(\mb{q}) \bar{\mb{h}}(\mb{\omega},\mb{\delta}, \mb{h}_{ga}, \mb{h}_{swr}), \label{Eq_PhysConstraints_mom}
\end{align}
\end{subequations}
where the satellite's inertial-frame angular momentum $\mb{h}_0 \in \R^3$ is conserved in the absence of external forces. Critically, \eqref{Eq_PhysConstraints} implicitly constrains $\mb{x}$ to a $(3m+3)$-submanifold $X(\mb{h}_0)$ of the ambient linear space $\R^{3m+7}$. Naturally, common dynamical approximations which compromise \eqref{Eq_PhysConstraints_mom} (e.g. modelling the MoI $J_{st}(\mb{\delta})$ as a constant) cannot continuously remain on $X(\mb{h}_0)$, ensuring that the corresponding solutions are not physical.

More critically, the state manifold $X(\mb{h}_0)$ implicitly constrains the local linear controllability of the CMG array: a feature which, if lost in approximation, can both slow trajectory optimization solvers and limit the practical efficacy of their solutions. 
To examine this, let $(\mb{z},\mb{v})$ capture local perturbations of $(\mb{x},\mb{u})$ in $\R^{3m+7}\times\R^{2m}$. The dynamics \eqref{Eq_Dynamics}, which satisfy \eqref{Eq_PhysConstraints} by construction, can then be expressed locally around $\mb{x}$ in the linear form
\begin{equation}
\label{Eq_LinDyn}
    \dot{\mb{z}} = A(\mb{x}) \mb{z} + B(\mb{x}) \mb{v},
\end{equation}
where $A(\mb{x}) \coloneqq \nicefrac{\partial f}{\partial x}$ and $B(\mb{x}) \coloneqq \nicefrac{\partial f}{\partial x}$. Similarly, a \supfirst-order Taylor expansion of \eqref{Eq_PhysConstraints} around $\mb{x}$ can be written in the compact form
$$
\begin{bmatrix}
    1 \\
    \bar{\mb{h}}(\mb{x}) \\
\end{bmatrix} \approx 
\begin{bmatrix}
    1 \\
    \bar{\mb{h}}(\mb{x}) \\
\end{bmatrix} +
\underbrace{
\begin{bmatrix}
    \mb{q}^\top / \norm{\mb{q}}   &    0  &   0    & 0 &  0   \\
    2 \left[ \bar{\mb{h}},  -\widehat{\bar{\mb{h}}} \right] O_L(\mb{q}^*)   &   A_s & J_{st} & D & A_g  \\
\end{bmatrix} 
}_{Z(\mb{x}) \coloneqq}
\mb{z}.
$$
For any $\mb{x}$, it can be shown that the rows of the matrix $Z(\mb{x})$ are linearly independent and span the null space of $A(\mb{x})$ (e.g. $A(\mb{x}) Z(\mb{x})^\top = \mb{0}$). That is, the local dynamics \eqref{Eq_LinDyn} are constrained to the tangent space $T_{\mb{x}} X$ of $X(\mb{h}_0)$ at $\mb{x}$ and can be rewritten in the reduced form
\begin{equation}
\label{Eq_ReducedLinDyn}
    \dot{\mb{s}} = \underbrace{ \left[ M(\mb{x}) A(\mb{x}) M(\mb{x})^\top \right] }_{A_s(\mb{x}) \coloneqq} \mb{s} + \underbrace{ \left[ M(\mb{x}) B(\mb{x}) \right]}_{B_s(\mb{x}) \coloneqq} \mb{v},    
\end{equation}
using the projected coordinate $\mb{s} \coloneqq M(\mb{x}) \mb{z} \in \R^{3m+3}$ where the rows of $M(\mb{x})$ are any orthonormal basis of $T_{\mb{x}} X$ (obtainable via \texttt{null}$(Z(\mb{x}))^\top$ in Matlab). Crucially, while the pair $(A_s,B_s)$ may be locally linearly controllable in the reduced $\R^{3m+3}$, \eqref{Eq_Dynamics} is \emph{not} locally linearly controllable in $\R^{3m+7}$. This feature is critical to the design of performant regulators, but is lost in many dynamical approximations.

\subsection{Singularities in Classic CMG Controllers}

Another common approximation of the satellite attitude dynamics \eqref{Eq_Dynamics} omits the CMG dynamics entirely, instead planning maneuvers using generalized body-torque commands $\mb{\tau}_r$ and the simplified attitude dynamics:
\begin{equation}
\label{Eq_TorqueDynamics}
\begin{aligned}
    \dot{\mb{q}} &=  \nicefrac{1}{2} \, O_L(\mb{q}) \, \widetilde{\mb{\omega}}\,, \\
    J\dot{\mb{\omega}} &= -\widehat{\mb{\omega}} \, J \mb{\omega} + \mb{\tau}_r \,,
\end{aligned}
\end{equation}
These pre-planned command torques are then converted to inputs for the CMG's using an appropriate Jacobian. For example, \cite{Oh1991} convert the commands $\mb{\tau}_r(t)$ to the minimum norm CMG gimbal rates $\dot{\mb{\delta}} = D^{\dagger} \mb{\tau}_r$ using the Moore-Penrose pseudo-inverse $D^{\dagger} \coloneqq D^{\top} (DD^\top)^{-1}$ of \eqref{Eq_Ddef}. Notably, this conversion uses \emph{only} the CMG gimbal motors to produce $\mb{\tau}_r$ as \eqref{Eq_TorqueDynamics} assumes a fixed (internally regulated) wheel speed. 

\begin{figure}
    \centering
    \includegraphics[width=3.2in]{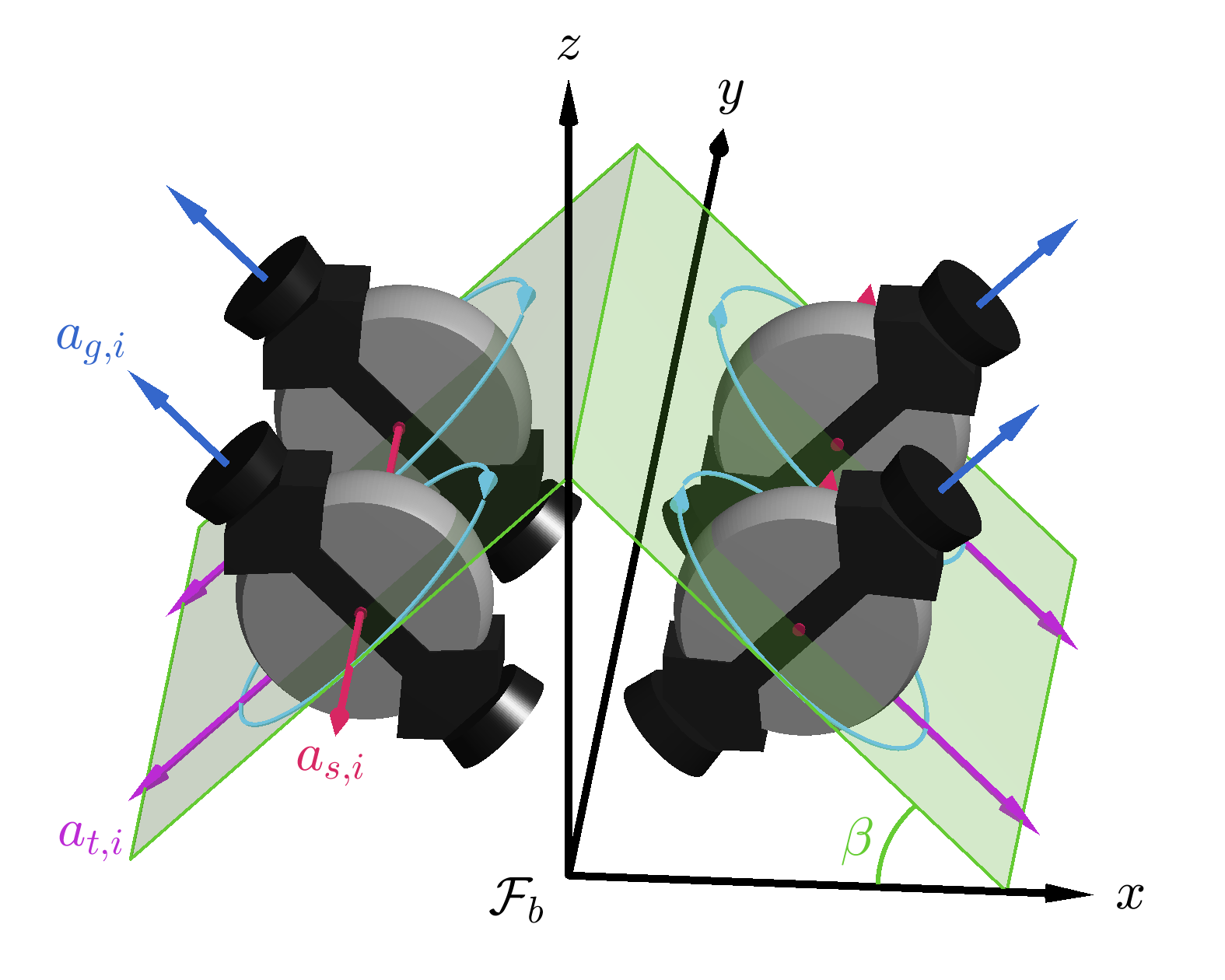}
    \caption{A 4-CMG array in rooftop configuration with an inclination of $\beta = \ang{45}$.
    \label{Fig_RooftopConfig}
    }
\end{figure}

\begin{figure}
    \centering
    \includegraphics[width=3.2in]{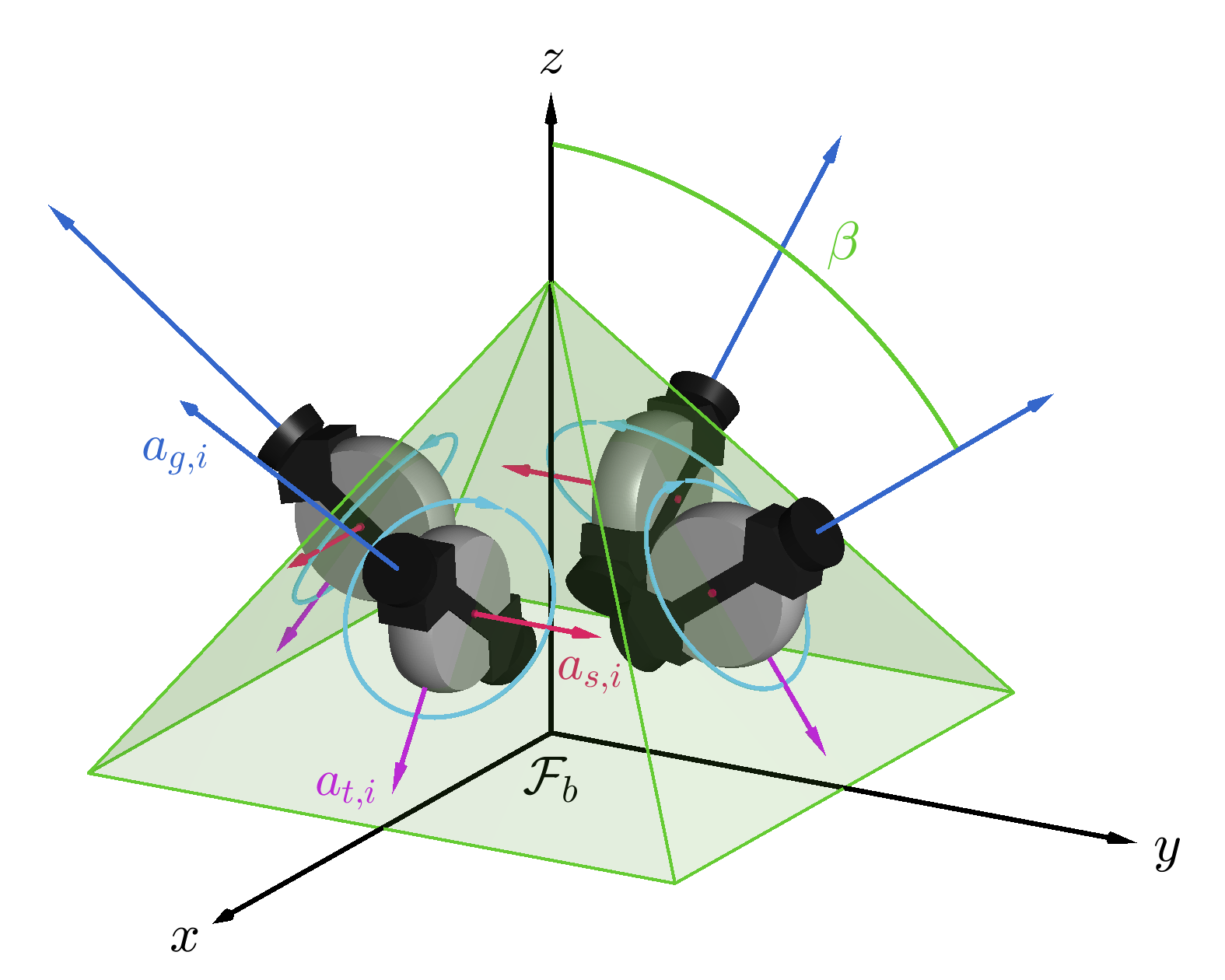}
    \caption{A 4-CMG array in Pyramid configuration with an inclination of $\beta = \ang{54.74}$.
    \label{Fig_PyramidConfig}
    }
\end{figure}

\begin{figure}
    \centering
    \includegraphics[width=2.0in]{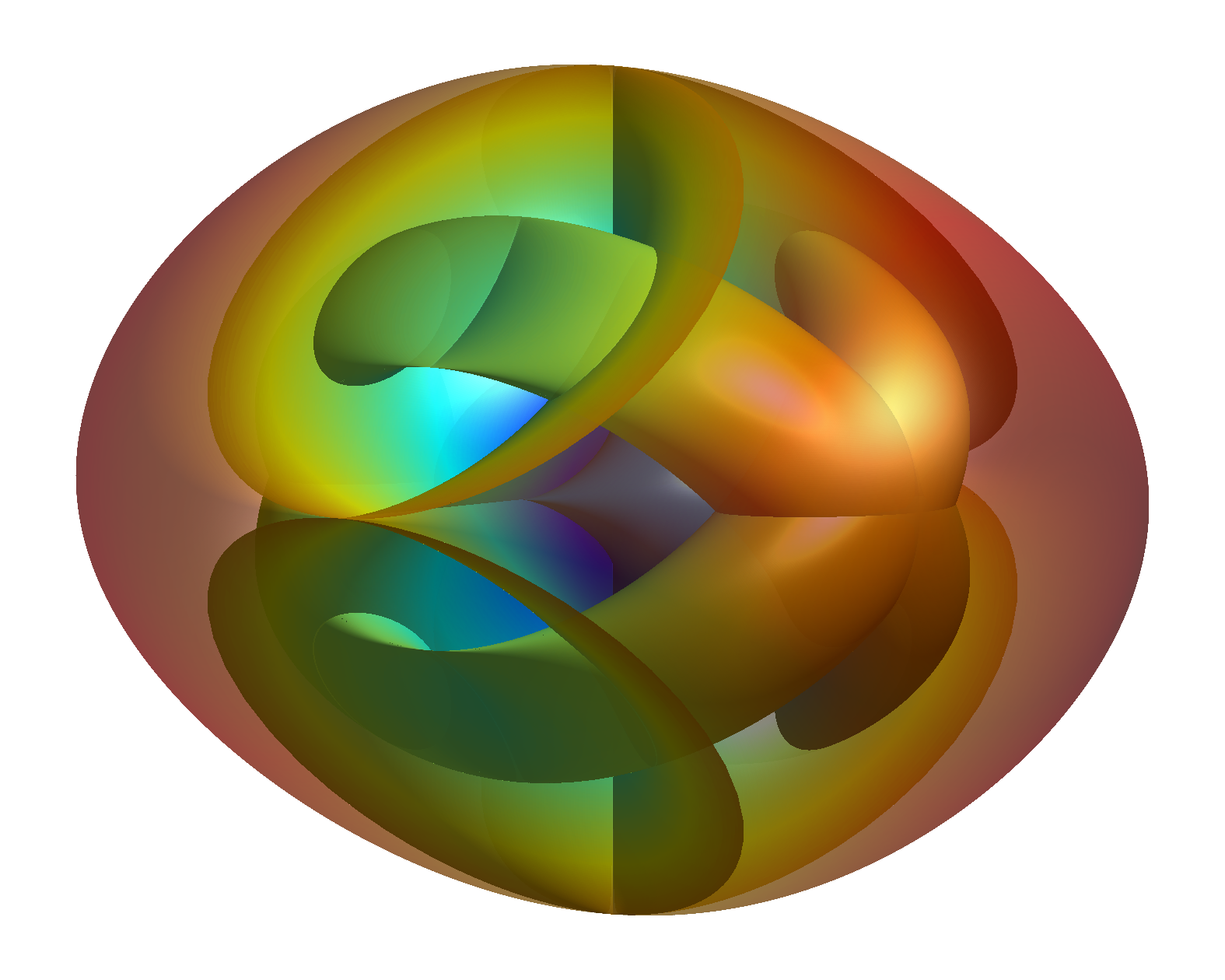}
    \caption{Singular configurations of a 6-CMG rooftop array projected into the array's momentum workspace.
    \label{Fig_RooftopSingSurf}
    }
\end{figure}

While this design assumption greatly simplifies the dynamics by separately regulating the CMG wheel speeds, the CMG gimbal motors alone are often insufficient to produce arbitrary command torques. Specifically, configurations $\mb{\delta}$ in which the matrix $A_t(\mb{\delta})$ is low rank \emph{cannot} produce torques within the missing torque space. For example, the default configurations (all $\delta_i=0$) of the rooftop and pyramid array geometries shown in Figures \ref{Fig_RooftopConfig} and \ref{Fig_PyramidConfig} cannot produce torques along the $y$ and $z$ axes respectively because all the CMG torque axes are coplanar. Such configurations are called \emph{singular} as the Jacobian $D$ shares rank with $A_s(\mb{\delta})$, producing a kinematic singularity in the above feedback strategy and making maneuvers in a neighborhood of these configurations inefficient or even impossible. Many approaches have been developed in the literature to avoid singularities, including popular strategies by \cite{Oh1991} and \cite{Schaub1998} which adaptively regularize the pseudo-inverse to avoid singularities or track a pre-computed set of `safe' configurations respectively. As shown by the projection in Fig. \ref{Fig_RooftopSingSurf}, these singular configurations are numerous in any array's configuration space and are highly dependent on the specifics of the array geometry. 


Naturally, these singularities pose significant challenges to the determination and integration of optimal control strategies. Primarily, strategies based only on \eqref{Eq_TorqueDynamics} can produce infeasible command torques for certain array configurations. Alternatively, formulations which incorporate the CMG control feedback to avoid singularities can slow or even stall solvers (due to the numerical ill-conditioning of $(DD^\top)^{-1}$ produced near singularities). Thankfully, VSCMG models like \eqref{Eq_Dynamics} can use the additional controllability provided by the CMG wheel motors to avoid the large majority of singular effects. Notably, kinematic singularities are \emph{always} present on the outer momentum envelope of the array shown in Fig. \ref{Fig_RooftopSingSurf} (where the CMG wheels saturate), but can be easily avoided with minor wheel speed regulation.

\section{Problem Formulation and Approach}

\subsection{Optimization Problem and Solver}
Having established physically accurate dynamics \eqref{Eq_Dynamics} free of the numerical effects of common dynamical approximations, we now introduce our optimal control problem. In this work, we examine the attitude transfer from $\mb{x}_0$ to $\mb{x}_d \in X(\mb{h}_0)$ while satisfying the dynamics $\dot{\mb{x}} = f(\mb{x},\mb{u})$ in \eqref{Eq_Dynamics}. For simplicity, we only consider rest-to-rest transfers between (non-singular) equilibrium points (e.g. $\mb{x}_0,\mb{x}_d \in \{x \in X(\mb{h}_0) : f(\mb{x},0) = 0,~\mb{\omega} = 0\}$). This general trajectory optimization problem is given by
\begin{equation}
\label{eq_Opt_Control_Prob}
\begin{aligned}
\min_{\mb{x}(\cdot),\mb{u}(\cdot)}  \quad&  \int_{0}^{T}{\ell(\mb{x}(t),\mb{u}(t)) \, dt} \;+\; m(\mb{x}(T)), \\
    \text{s.t.} \quad &\dot{\mb{x}} = f(\mb{x},\mb{u}), \; \mb{x}(0) = \mb{x}_0,
\end{aligned}
\end{equation}
where the notation $\mb{x}(\cdot)$ denotes the entire curve $\mb{x}(t)$ over the interval $t\in[0,T]$. To choose the stage and terminal cost functionals $\ell$ and $m$ and solve \eqref{eq_Opt_Control_Prob}, we adapt the approach in \cite{Dearing2021} using the \textit{PR}ojection-\textit{O}perator-based \textit{N}ewton's method for \textit{T}rajectory \textit{O}ptimization (PRONTO). Briefly, PRONTO is a direct method based on a modified Newton descent step using \supsecond-order approximations of the local dynamics and cost function. In particular, PRONTO's continuous-time solution iterates (generated via numerical integration) implicitly satisfy the system's dynamics, enabling more effective descent directions to be obtained from a significantly reduced search space. This aspect makes PRONTO particularly effective in problems with highly complex dynamics. We refer the reader to \cite{Dearing2021} for an in-depth discussion of PRONTO and the specifics regarding its implementation.


\subsection{Cost Co-Design for PRONTO Solver}

The constrained state manifold $X(\mb{h}_0)$ generated by \eqref{Eq_PhysConstraints} presents a challenge when choosing the stage and terminal cost functionals $\ell$ and $m$. Following the strategies in \cite{Dearing2021}, a simple but effective choice for these functions are the quadratic forms 
\begin{subequations}
\label{eq_CostDefs}
\begin{align}
    \ell(\mb{x},\mb{u}) &\coloneqq \tfrac{1}{2} \norm{\mb{x} - \mb{x}_d}^2_{Q(\mb{x}_d)} + \tfrac{1}{2} \norm{\mb{u}}^2_R, \label{eq_CostDefs_l}\\
           m(\mb{x}) &\coloneqq \tfrac{1}{2} \norm{\mb{x} - \mb{x}_d}^2_{P(\mb{x}_d)}, \label{eq_CostDefs_m}
\end{align}
\end{subequations}
where $\norm{\mb{v}}^2_R$ is shorthand for the semi-norm $\mb{v}^\top R \,\mb{v}$. Specifically, suitable positive semi-definite quadratic state weights $Q(\mb{x}_d),P(\mb{x}_d) \in \R^{(3m+7)\times(3m+7)}$ and positive definite control weight $R \in \R^{2m\times2m}$ are produced in the design of a locally exponentially stabilizing Linear Quadratic Regulator (LQR) around the target state $\mb{x}_d$. Recalling that the dynamics \eqref{Eq_Dynamics} are implicitly constrained to the state manifold $X(\mb{h}_0)$ (and not locally linearly controllable), this regulator must be designed for the reduced dynamics \eqref{Eq_ReducedLinDyn} on the tangent space $T_{\mb{x}_d} X$ at $\mb{x}_d$. Specifically, given symmetric positive definite matrices $Q_s\in \R^{(3m+3)\times (3m+3)}$ and $R \in \R^{2m\times2m}$, there will be a positive definite matrix $P_s \in \R^{(3m+3)\times (3m+3)}$ satisfying the Algebraic Riccati Equation (ARE): 
\begin{equation}
\label{Eq_ARE}
    Q_s = P_s B_s R^{-1} B_s^\top P_s - A_s^\top P_s - P_s A_s.
\end{equation}
We can then lift the state cost $Q_s$, the Riccati solution $P_s$, and (if desired) the feedback regulator $K_s = R^{-1} B_s^\top P_s$ into the ambient space using the projection $M(\mb{x}_d)$:
\begin{equation}
\label{Eq_QPK_lift}
\begin{aligned}
    Q(\mb{x}_d) &= M(\mb{x}_d)^\top Q_s M(\mb{x}_d), \\
    P(\mb{x}_d) &= M(\mb{x}_d)^\top P_s M(\mb{x}_d), \\
    K(\mb{x}_d) &= K_s M(\mb{x}_d),
\end{aligned}
\end{equation}
where, by design, $Q(\mb{x}_d)$ and $P(\mb{x}_d)$ are positive definite on $T_{\mb{x}_d} X$ and zero on its orthogonal complement. Note that this design strategy requires the pair $(A_s(\mb{x}_d), B_s(\mb{x}_d))$ to be linearly controllable at the target $\mb{x}_d$. While this is certainly true for the problems discussed below, a general proof for arbitrary array geometries and target states $\mb{x}_d \in X(\mb{h}_0)$ under the dynamics \eqref{Eq_Dynamics} is nontrivial and absent from existing literature. However, \cite{Bhat2015} have shown a comparatively underactuated variation of \eqref{Eq_Dynamics} to be linearly controllable around non-singular equilibrium points. This, combined with our own extensive numerical testing, indicate that this restriction is unlikely to have any significant practical impact.



\section{Numerical Evaluation}
In this section, we examine the mean statistics for PRONTO solutions to the optimization problem \eqref{eq_Opt_Control_Prob} as well as specific solution features for an additional model problem. In particular, we examine the popular rooftop and pyramid array geometries shown in Figures \ref{Fig_RooftopConfig} and \ref{Fig_PyramidConfig} respectively with the platform inertias (in \SI{}{\kilo\gram \meter^2})
$$
\begin{aligned}
        \diag(J) &= 
        \begin{bmatrix}
            1500   &    1500  &   2000  
        \end{bmatrix},
    &
        J_g &= 0.115,
    \\
        J_{sw} &= 0.075, \, J_{sg} = 0.015,
    &
        J_t &= 0.001,
\end{aligned}
$$
and a target (and initial) CMG wheel momentum of $h_{swr,t} = \SI{25}{\kilo\gram \meter^2 \second^{-1}}$. For the generation of our LQR cost functional (and the projection regulator used by the PRONTO solver) following \eqref{Eq_ARE} and \eqref{Eq_QPK_lift}, the positive definite weight matrix $Q_s \coloneqq M(\mb{x}_d) Q_c M(\mb{x}_d)^{\top}$ on the controllable subspace can be generated using the structure
$$
Q_c \coloneqq \mathrm{diag}([ \rho_{q} \mb{1}_4 \,;\, \rho_{hswr} \mb{1}_m \,;\, \rho_{\omega} \mb{1}_3 \,;\, \rho_{\delta} \mb{1}_m \,;\, \rho_{hga}  \mb{1}_m  ]),
$$
to manage individual state error weights via the scalar weights $\rho_i$, where $\mb{1}_m\in\R^m$ denotes a vector of 1's. For the cost functional and regulator respectively, these weights were chosen to be
$$
\begin{aligned}
    [\rho_q, \rho_{hswr}, \rho_{\omega}, \rho_{\delta}, \rho_{hga}]_{\mathrm{cost}} &= [ 5, 10, 0.1, 0.01, 50 ], \\
    [\rho_q, \rho_{hswr}, \rho_{\omega}, \rho_{\delta}, \rho_{hga}]_{\mathrm{reg}} &= [ 3\cdot10^{4}, 3, 200, 0.3, 3]\cdot10^{-4}. 
\end{aligned}
$$
The control weights for $R \coloneqq \diag([\rho_{ug} \mb{1}_m; \rho_{uw} \mb{1}_m])$ for the cost function and regulator were likewise chosen as
$$
[\rho_{ug}, \rho_{uw}]_{\mathrm{cost}} = [1,1], \qquad [\rho_{ug}, \rho_{uw}]_{\mathrm{reg}} = [1,3]\cdot10^{-5}.
$$
Finally, initial guess trajectories provided to PRONTO were generated using the well-known \emph{Singularity Robust} feedback control law presented in \cite{Oh1991}, with a time horizon of \SI{180}{\second} found to allow sufficient convergence for single rest-to-rest rotations. All solutions for both solvers were computed in Matlab on an AMD 5800X CPU platform with 32 GBs of 32 MHz memory.


\subsection{Mean Performance comparison}
In order to compare the performance of the initial feedback solution from \cite{Oh1991} to the optimal trajectory determined by PRONTO, solutions were computed for 10 randomly generated rest-to-rest attitude transfers, with $\mb{q}_0,\mb{q}_d \in S^3$, $\mb{\omega}_0, \mb{\omega}_d = \mb{0}$, and the remaining momentum states chosen as non-singular zero momentum ($\mb{h}_0 = 0$) configurations for that array satisfying $\mb{x}_0, \mb{x}_d \in X(\mb{h}_0)$. Mean performance statistics comparing the initial guess and optimizer over these runs are shown in Table \ref{Tab_GuessOptStats}.

\begin{table}
\label{Tab_GuessOptStats}
\caption{Optimal Trajectory Statistics}
\centering
\begin{tabular}{|l|rr||rr|}
    \hline
                                                                         & \multicolumn{2}{c||}{Rooftop}          & \multicolumn{2}{c|}{Pyramid}          \\ \hline
    \multicolumn{1}{|c|}{Metric}                                         & \multicolumn{1}{c|}{Guess}  & \multicolumn{1}{c||}{Opt.} & \multicolumn{1}{c|}{Guess}    & \multicolumn{1}{c|}{Opt.} \\ \hline
    Comp. Time $[\SI{}{\minute}]$                                        & \multicolumn{1}{r|}{NA}      & 15.45                        & \multicolumn{1}{r|}{NA}        & 36.47                        \\ \hline
    Maneuver Cost                                               & \multicolumn{1}{r|}{83.56}  & 39.90                        & \multicolumn{1}{r|}{91.53}    & 34.07                        \\ \hline
    Control Effort $[\SI{}{\newton \meter \second}]$                     & \multicolumn{1}{r|}{109.77} & 25.40                        & \multicolumn{1}{r|}{140.21}   & 26.05                        \\ \hline
    Maneuver Energy $[\SI{}{\joule}]$                                    & \multicolumn{1}{r|}{5.07}   & 15.15                        & \multicolumn{1}{r|}{4.82}     & 19.95                        \\ \hline
    Maneuver Time $[\SI{}{\second}]$                                     & \multicolumn{1}{r|}{95.70}  & 47.93                        & \multicolumn{1}{r|}{100.35}   & 37.39                        \\ \hline
    Final Att. Error $[\SI{}{\degree}]$                              & \multicolumn{1}{r|}{0.83}   & 0.06                         & \multicolumn{1}{r|}{1.68}     & 0.10                         \\ \hline
    Max $\mb{u}_g$ $[\SI{}{\newton \meter}]$                             & \multicolumn{1}{r|}{0.47}   & 1.19                         & \multicolumn{1}{r|}{1.5E-4} & 3.9E-3                     \\ \hline
    Max $\mb{u}_w$ $[\SI{}{\newton \meter}]$                             & \multicolumn{1}{r|}{0.50}   & 1.08                         & \multicolumn{1}{r|}{1.2E-4} & 4.0E-3                     \\ \hline
\end{tabular}
\end{table}

Examining Table \ref{Tab_GuessOptStats}, we first note that solutions to this challenging problem are not obtained easily. For the rooftop and pyramid geometries, PRONTO takes an average of 15 and 35 minutes respectively to reduce the objective cost to 47\% and 42\% of its original value. Notably, Matlab's \texttt{ode45} function limits the algorithm to a single CPU thread (a limitation shared by the majority of spacecraft CPU's). Interestingly, solutions for the pyramid geometry were far more computationally expensive than those of the rooftop array, indicating a higher intrinsic complexity in the effective operation of that geometry.

While the maneuver convergence time and terminal attitude error show similar reductions to that of the objective function, a more interesting effect is observed in the mean maneuver efficiency. In particular, while the total control effort (integration of $\sum_i\abs{u_i}$) shows a reduction of 77\% and 81\% respectively, the optimal maneuver uses \emph{far} more true electric power (3-4x) than the initial guess. While this increased energy usage partially originates from an aggressive cost function weighting, it also highlights a shortcoming of the cost functional \eqref{eq_Opt_Control_Prob} for this system. Specifically, recall that the power $P = \tau_w \omega_w$ used by an electric motor is proportional to its shaft speed. Since both the CMG gimbal and wheel have \emph{variable} speeds in practice, the cost penalty on $\abs{u}$ does not directly penalize the array's true power consumption. To the author's knowledge, this subtlety has not been addressed in existing literature and suggests substantial performance improvements.


\subsection{Optimal Trajectory Features}
\begin{figure}
    \centering
    \includegraphics[width=3.2in]{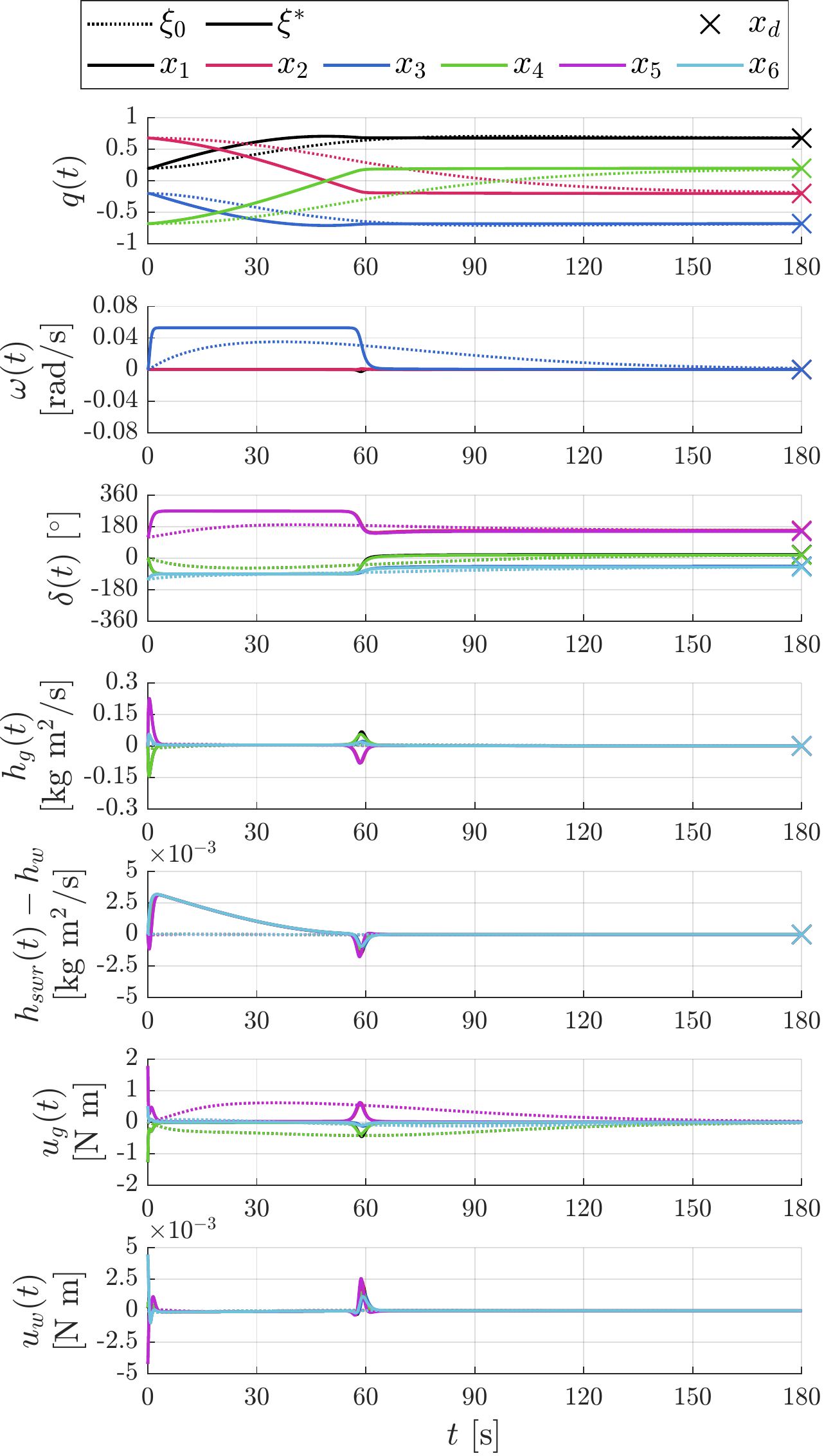}
    \caption{Guess ($\xi_0$) and Optimal ($\xi^*$) Trajectories for a \ang{180} z-axis rotation of the Rooftop CMG geometry.
    \label{Fig_Rooftop_Model_Trajectory}
    }
\end{figure}

\begin{figure}
    \centering
    \includegraphics[width=3.2in]{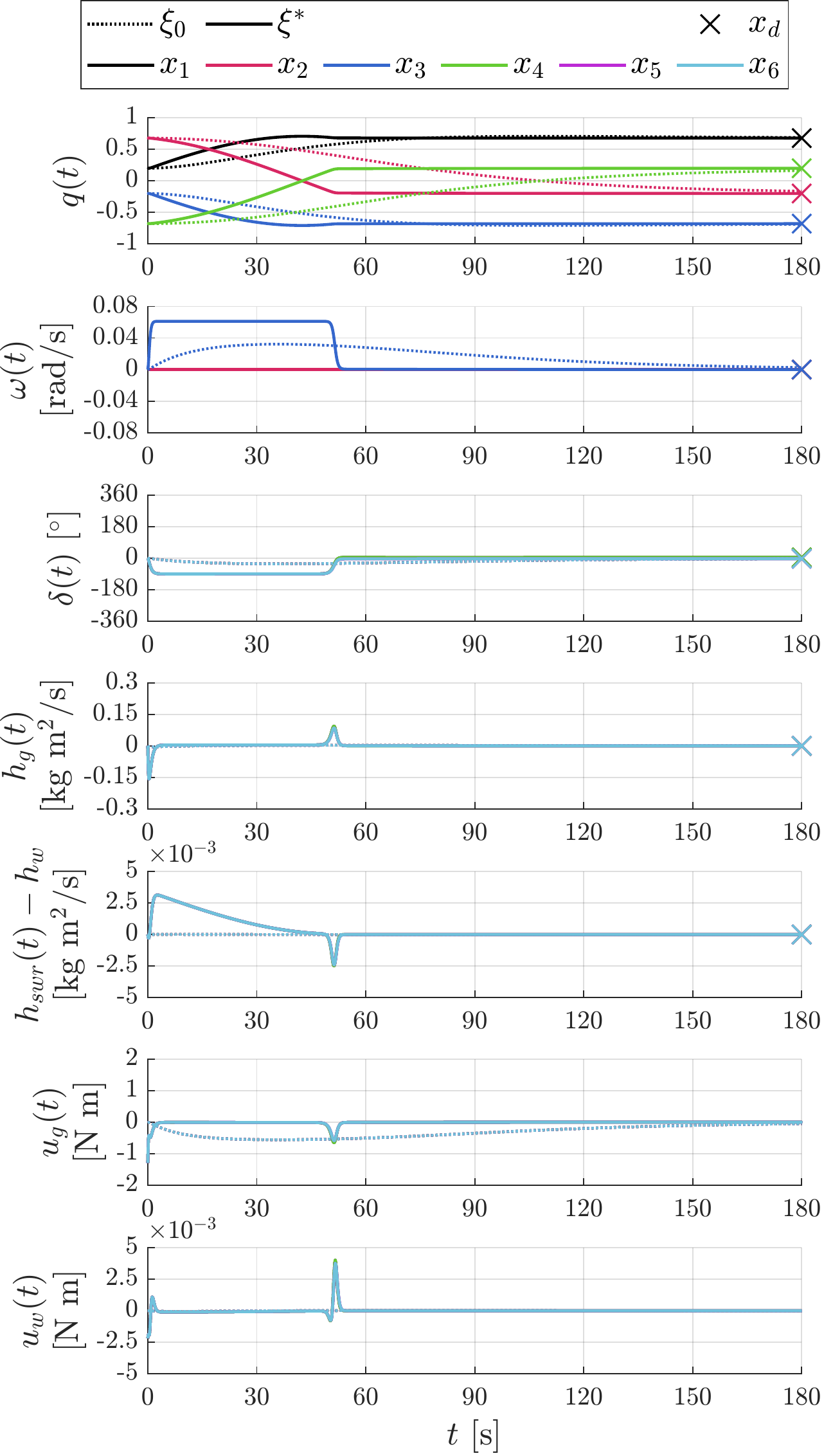}
    \caption{Guess ($\xi_0$) and Optimal ($\xi^*$) Trajectories for a \ang{180} z-axis rotation of the Pyramid CMG geometry.
    \label{Fig_Pyramid_Model_Trajectory}
    }
\end{figure}

Examining these solutions in more detail, Figures \ref{Fig_Rooftop_Model_Trajectory} and \ref{Fig_Pyramid_Model_Trajectory} show the optimal and guess trajectories for a \ang{180} rotation about the $z$-axis for the rooftop and pyramid geometries respectively. These optimal trajectories display several interesting features
\begin{enumerate}
    \item $\mb{\omega}(t)$ saturates in both maneuvers and geometries.
    \item Both $\mb{u}_w$ and $\mb{u}_g$ are impulsive in nature.
    \item The CMG angles $\delta_i$ display unusual coordination.
\end{enumerate}
Regarding observation 1, we remind the reader that slew rate constraints are \emph{not} considered in this problem. Instead, the apparent maximum rotation rate results from the finite momentum capacity of a CMG array (the envelope of Figure \ref{Fig_RooftopSingSurf}), which enforces a maximum rotation rate along any axis. This natural property of MED's also informs upon observation 2, with each control input acting impulsively to rapidly achieve the array configuration for this maximum rotation rate. Regarding observation 3, the coordination of the pyramid array is intuitive as its symmetry with the requested rotation axis clearly promote symmetry in the actuators. However, the coordination for the rooftop geometry is far more interesting. While we might expect the CMG's to coordinate in groups with shared gimbal axes (as sides of the rooftop), they instead operate in pairs \emph{across} the rooftop. This intriguing behavior was observed for multiple maneuvers with different rotation axes and warrants further investigation.


\section{Conclusions}
In this work, we developed a numerically tractable trajectory optimization problem for rest-to-rest attitude transfers with CMG-driven spacecraft. This included the development of a specialized dynamical model which, while more complex and nonlinear than traditional approximated models, avoids many of the practical complications which slow or stall conventional solvers. To develop and solve this specialized trajectory optimization problem, we designed a locally stabilizing LQR on the system's configuration manifold, then lifted it into the ambient state space to produce suitable terminal and running LQ cost  functionals. Finally, we examined the performance benefits and drawbacks of solutions to this optimization problem: an investigation which revealed both significant performance improvements under our formulation, potential avenues for future performance improvements, and interesting solution features which could inform the development of future control laws.

\bibliography{References}       

\end{document}